\documentclass[12pt]{article}
\usepackage{amssymb,amsmath,amsfonts,t1enc}
\title{On Combinatorial Problem Concerning Partitions of a~Box into Boxes}
\author{Apoloniusz Tyszka}
\date{}
\begin{document}
\def\thefootnote{}
\maketitle
\begin{abstract}
We consider partitions of $n$-dimensional boxes in $\mathbb{R}^n$, $n \geq 2$,
into a~finite number of boxes with pairwise disjoint interiors. We study sets
$X \subseteq (0,\infty)$ with the Property $(W_n): $ for every $n$-dimensional
box $P$ and every partition of $P$, if each constituent box has one side with
the length belonging to $X$, then the length of one side of $P$ belongs to~$X$.
We prove that the set $X\subseteq (0,\infty)$ has Property $(W_n)$ if and only 
if $X$ is closed with respect to the operations: $x+y$ and $x+y+z-2\min(x,y,z)$.
\end{abstract}
\par
We consider partitions of $n$-dimensional boxes in $\mathbb{R}^n$, $n \geq 2$,
into a~finite number of $n$-dimensional boxes with pairwise disjoint interiors.
We study sets $X\subseteq(0,\infty)$ with the Property $(W_n):$ for every
$n$-dimensional box $P$ and every partition of $P$, if each constituent box has
one side with the length belonging to $X$, then the length of one side of
$P$ belongs to~$X$.
\vskip 0.2truecm
\par
\noindent
\textbf{Theorem\ }(see: \cite{5} for groups, \cite{4} for generalizations
of groups, \cite{1}--\cite{3} for earlier results).\ \emph{The set $X \subseteq (0,\infty)$ has Property $(W_n)$ if and only if $X$
is closed with respect to the operations: $x+y$ and $x+y+z-2\min(x,y,z)$.}
\vskip 0.2truecm
\par
\noindent
\textit{Proof.} Necessity: We can assume that $x>z,\ y>z$. Case $n=2$:
the result follows from the following partitions of a square:
\begin{center}
\begin{picture}(80,35)(0,-10)
\put(0,0){\line(0,1){20}}
\put(0,20){\line(1,0){20}}
\put(20,20){\line(0,-1){20}}
\put(20,0){\line(-1,0){20}}
\put(15,0){\line(0,1){20}}
\put(5,-3){$x$}
\put(17,-3){$y$}
\put(22,10){$x+y$}
\put(60,0){\line(0,1){20}}
\put(60,20){\line(1,0){20}}
\put(80,20){\line(0,-1){20}}
\put(80,0){\line(-1,0){20}}
\put(70,-3){$x$}
\put(57,7){$x$}
\put(65,22){$y$}
\put(82,15){$y$}
\put(70,20){\line(0,-1){10}}
\put(72,12){$z$}
\put(63,17){\line(0,-1){17}}
\put(63,10){\line(1,0){17}}
\put(66,7){$z$}
\put(60,17){\line(1,0){10}}
\end{picture}
\end{center}
Case $n>2$: we multiply the first square by $[0, x + y]^{n-1}$, the second
square by $[0,x+ y - z]^{n-1}$. Sufficiency: we first prove (see Eulerian path
method in \cite{5} by Michael S. Paterson, Univ. of Warwick, Coventry, England)
that if each constituent box from the partition of the box $P$ has one side with
the length belonging to $X$, then there exist points $Y_1, Y_2,\ldots, Y_m$
lying on the one side of $P$ such that this side is equal to the segment
$\overline{Y_1Y_m}$ and every distance $|Y_iY_{i+1}|$ belongs to $X$. We choose
the cartesian coordinate system with axes $x_j (j = 1, 2, \ldots, n)$ which
are parallel to the sides of the box $P$. For each constituent box $P_k$ we
choose $c(k) \in \{1, 2, \ldots, n\}$ such that the length of $P_k$ in direction
$x_{c(k)}$ belongs to $X$. We define an undirected graph $G$ in the following
way: as a~vertex set we put the set of all vertices of constituent boxes, as
an edge set we put the set of all pairs $(s, P_k)$, where $s$ is a side of
$P_k$ lying in the direction $x_{c(k)}$. Each vertex of $G$ (except the vertices
of the box $P$ which lie on $1$ edge) is a vertex of an even number of
constituent boxes, hence it lies on an even number of edges of $G$. Thus a~walk
away along edges that begins at one vertex of $P$ and does not repeat any edges
will not terminate until it hits another vertex of $P$. This observation produce points $Y_1,Y_2,\ldots,Y_m$.
\vskip 0.2truecm
\par
Now we are ready to prove that $|Y_1Y_m| \in X$. It suffices to prove that if
there exist $m$ points with the required property (for $m > 3$), then there
exist $m - 1$ or $m - 2$ points with this property. We can assume
$Y_1,Y_2,\ldots Y_m$ are different and for every $i$ such that $1<i<m$,
$Y_i\not\in\overline{Y_{i-1}Y_{i+1}}$, hence $Y_3\in\overline{Y_1Y_2}$ and
$Y_{m-2}\in\overline{Y_{m-1}Y_m}$, so $|Y_2Y_3|<|Y_1Y_2|$ and
$|Y_{m-1}Y_m|>|Y_{m-2}Y_{m-1}|$. We choose the smallest $i>2$ such that
$|Y_iY_{i+1}|>|Y_{i-1}Y_i|$. From this choice of $i$ we have
$|Y_{i-1}Y_i|<|Y_{i-2}Y_{i-1}|$, so
\begin{center}
\begin{picture}(83,10)(0,-3)
\linethickness{0.5pt}
\put(-2,0){\circle*{1}}
\put(-1,2){$Y_{i-2}$}
\put(0,0){\line(1,0){30}}
\put(32,0){\circle*{1}}
\put(33,2){$Y_i$}
\put(34,0){\line(1,0){10}}
\put(46,0){\circle*{1}}
\put(47,2){$Y_{i-1}$}
\put(48,0){\line(1,0){30}}
\put(80,0){\circle*{1}}
\put(81,2){$Y_{i+1}$}
\end{picture}
\end{center}
$|Y_{i-2}Y_{i+1}|=|Y_{i-2}Y_{i-1}|-|Y_{i-1}Y_{i}|+|Y_{i}Y_{i+1}|=$
\\
$\makebox[3em]{}|Y_{i-2}Y_{i-1}|+|Y_{i-1}Y_{i}|+|Y_{i}Y_{i+1}|-2\min(|Y_{i-2}Y_{i-1}|,|Y_{i-1}Y_{i}|,|Y_{i}Y_{i+1}|).$
\\
\par
\noindent
$X$ is closed with respect to the operation $x+y+z-2\min(x,y,z)$, hence
$|Y_{i-2}Y_{i+1}|\in X$ and we can reduce the number of points to $m-2$, this
ends the proof of sufficiency.
\vskip 0.2truecm
\par
Our proof is now complete.
\newline
\rightline{$\Box$}

Apoloniusz Tyszka\\
Technical Faculty\\
Hugo Ko\l{}\l{}\k{a}taj University\\
Balicka 116B, 30-149 Krak\'ow, Poland\\
E-mail address: {\it rttyszka@cyf-kr.edu.pl}
\end{document}